\RequirePackage{fix-cm}
\documentclass[smallextended]{svjour3}
\smartqed
\usepackage[tracking=true]{microtype}
\usepackage[utf8]{inputenx}
\usepackage{mathptmx}
\usepackage{array}                                                              
\usepackage{booktabs}
\usepackage{amsmath}
\usepackage{amssymb}
\usepackage{mathtools}
\usepackage[oldcommands,algoruled,longend]{algorithm2e}
\usepackage[sort,compress]{cite}
\usepackage{graphicx}
\SetKwInput{KwInit}{Initialize}

\DeclareMathOperator*{\argmin}{arg\,min}

\DeclareMathOperator{\gshrink}{gshrink}
\DeclareMathOperator{\shrink}{shrink}
\providecommand{\abs}[1]{\left|#1\right|}
\providecommand{\norm}[1]{\left\|#1\right\|}
\providecommand{\set}[1]{\left\lbrace #1\right\rbrace}
\providecommand{\coloneqq}{\mathrel{\mathop:}=}
\journalname{}
\begin{document}
\title{%
  Bregman Iteration for Correspondence Problems
}
\subtitle{A Study of Optical Flow}
\author{Laurent Hoeltgen\and{}Michael Breuß}
\institute{%
  Chair for Applied Mathematics,
  Brandenburg Technical University, Cottbus-Senftenberg, Germany\\
  \email{\ensuremath{\lbrace{}}hoeltgen, breuss\ensuremath{\rbrace{}}@b-tu.de}%
}
\date{\today}
\maketitle
\begin{abstract}
      Bregman iterations are known to yield excellent results for denoising,
      deblurring and compressed sensing tasks, but so far this technique has
      rarely been used for other image processing problems. In this paper we
      give a thorough description of the Bregman iteration, unifying thereby
      results of different authors within a common framework. Then we show how
      to adapt the split Bregman iteration, originally developed by Goldstein
      and Osher for image restoration purposes, to optical flow which is a
      fundamental correspondence problem in computer vision. We consider some
      classic and modern optical flow models and present detailed algorithms
      that exhibit the benefits of the Bregman iteration. By making use of the
      results of the Bregman framework, we address the issues of convergence and
      error estimation for the algorithms. Numerical examples complement the
      theoretical part.
      \keywords{Bregman Iteration\and Split Bregman Method\and Optical Flow}
      \subclass{MSC~65Kxx, MSC~65Nxx}
\end{abstract}
\section{Introduction}
\label{sec:1}
%
In 2005, Osher~\emph{et al.}~\cite{Osher2005} proposed an algorithm for the
iterative regularisation of inverse problems that was based on findings of
Bregman~\cite{Bregman1967}. They used this algorithm, nowadays called Bregman
iteration, for image restoration purposes such as denoising and deblurring.
Especially in combination with the Rudin-Osher-Fatemi (ROF) model for denoising
\cite{ROF92} they were able to produce excellent results. Their findings caused
a subsequent surge of interest in the Bregman iteration. Among the numerous
application fields, it has for example been used to solve the basis pursuit
problem \cite{Cai2009,Osher2008,Yin2007} and was later applied to medical
imaging problems in \cite{Lin2006}. Further applications include deconvolution
and sparse reconstructions \cite{ZBBO09}, wavelet based denoising \cite{Xu2006},
and nonlinear inverse scale space methods \cite{Burger2006,Burger2005}. An
important adaptation of the Bregman iteration is the split Bregman method (SBM)
\cite{Goldstein2009} and the linearised Bregman approach \cite{Cai2009}. The SBM
can be used to solve $L_1$-regularised inverse problems in an efficient way. Its
benefits stem from the fact that differentiability is not a necessary
requirement on the underlying model and that it decomposes the original
optimisation task in a series of significantly easier problems that can be
solved very efficiently, especially on parallel architectures. The Bregman
algorithms belong to the family of splitting schemes as well as to the
primal-dual algorithms \cite{E2009,Setzer2009} which enjoy great popularity in
the domain of image processing and which are still a very active field of
ongoing research \cite{ODBP2015,GLY2015}.\par
The aim of this paper is to contribute to the mathematical foundation of the
rapidly evolving area of computer vision. We explore the use of the Bregman
framework, especially the application of the split Bregman method, for the
problem of optical flow (OF) which is of fundamental importance in that field,
cf.\ \cite{Aubert2006,KSK98,TV98}. We give a thorough discussion of the Bregman
framework, thereby unifying results of several recent works. Then we show how to
adapt the SBM to several classic and modern OF models. Detailed descriptions of
corresponding algorithms are presented. Employing the Bregman framework, we show
that convergence for these methods can be established and error estimates can be
given.\par
%
\subsection{The optical flow problem}
\label{sec:11}
%
The OF problem is an ill-posed inverse problem. It consists in determining the
displacement field between different frames of a given image sequence by looking
for correspondences between pixels. In many cases such correspondences are not
unique or simply fail to exist because of various problems such as noise,
illumination changes and overlapping objects. Nevertheless, the study of the OF
problem is of fundamental importance for dealing with correspondence problems
such as stereo vision where accurate flow fields are necessary
\cite{BS07,MMK98,SBW2005,MJBB2015}. For solving the OF problem in a robust way,
variational formulations and regularisation strategies belong to the most
successful techniques. Those methods have been studied for almost three decades,
starting from the approach of Horn and Schunck \cite{HS81}. During this period
of time, many efforts have been spent to improve the models cf.\
\cite{BA96,ZBWS09,BBPW04,MP02,NBK08,WCPB09,WPB10,XJM10,ZBWVSRS09,BZW2011,ZBW2011}
for an account of that field.\par{}
While many developments have been made on the modelling side, there are just a
few works concerned with the mathematical validation of algorithms. In
\cite{ki08,mm04} it has been shown that the classic numerical approach of Horn
and Schunck converges. Furthermore, the authors of \cite{mm04} showed that the
linear system obtained through the Euler-Lagrange equations has a symmetric and
positive definite matrix and thus allows the usage of many efficient solvers.
The authors of \cite{Wedel2008,Pock2007} developed an algorithm that solves the
so called TV-$L^{1}$ model through an alternating minimisation scheme. This is
applied to a variational formulation that augments the original energy
functional with an additional quadratic term. This quadratic term allows the
authors to divide their objective into simpler subproblems for which efficient
solvers exist. In practice their approach yields excellent results. However, in
general it does not converge towards the solution of the original energy
functional but to a solution of the augmented variational formulation.
Alternative approaches to minimise the occurring variational models include
\cite{CP2011,OCBP2014,OBP2015}. These well performing algorithms possess good
convergence properties, but may require additional regularity conditions, such
as strong convexity of the considered energy, see \cite{OBP2015}. The author of
\cite{br06} discusses the usage of efficient algorithms such as the Multigrid
approach \cite{B1973,BL2011,BHM00,Hac85,Wes92} and the so called
Lagged-Diffusivity or Ka\u{c}anov method \cite{CM99,FKN73,KNPS68}. Finally, it
is also possible to consider the solutions of the Euler-Lagrange equations as a
steady-state of a corresponding diffusion-reaction system that one may solve by
means of a steepest descend approach \cite{wb05,WBPB04}. Recent developments
have extended the study of the OF problem onto dynamic non-Euclidean settings
where the motion is estimated on an evolving surface \cite{KLS2015}. Other
trends include the use of powerful (deep) learning strategies
\cite{SRLB2008,DFIH2015} or even combinations of variational and machine
learning approaches \cite{WRHS2013}.\par
%
\subsection{Our contribution}
\label{sec:12}
%
We present an approach to the OF problem by exploring the Bregman framework.
Despite their usefulness, Bregman iterations have received little attention in
the context of OF up to now. Early attempts inlcude \cite{LBS2010,H2010}. Here,
we propose mathematically validated methods for OF models, among them the
prominent model of Brox \emph{et al.}~\cite{BBPW04}.\par
The main contribution of this work lies in the thorough presentation of the
Bregman framework and the proof of convergence of the algorithms in the context
of OF, thus giving the numerical solution of the OF problem a solid mathematical
basis. To this end, we adapt the general convergence theory of the Bregman
framework to the OF algorithms and show that the SBM iterates converge towards a
minimiser of the considered energy functionals. Related questions that are
important in the context of the numerical processing will also be discussed
here. For instance, we will show that the arising linear systems have a
symmetric and positive definite matrix. The assumptions for this are quite weak
and naturally met in almost all setups.\par
%
\subsection{Paper Organisation}
\label{sec:13}
%
In Section~\ref{sec:2}, we give a brief account of mathematical prerequisites,
whereas Section~\ref{sec:3} elaborates on the Bregman framework. Next, in
Section~\ref{sec:4} we give an account of the OF models we consider, and how to
formulate the corresponding algorithms in terms of the SBM. Finally, we
complement the theoretical developments by some numerical experiments given in
Section~\ref{sec:5} and finish the paper by some conclusions.
%
\section{Mathematical prerequisites}
\label{sec:2}
%
In this work we strongly rely on the notion of subdifferentiability as it grants
us the ability to handle non-differentiable robust regularisers, such as the
$\ell_{1}$ norm, in similar style as smooth functions. For a thorough analysis
of this important concept in convex optimisation we refer to the excellent
presentations in \cite{Ekeland1999,Rockafellar1997,rj98}. Here, we merely recall
the definition of a subdifferential. The subdifferential
$\partial\varphi\left(\bar{x}\right)$ of a function
$\varphi\colon{}\mathbb{R}^{n}\to\mathbb{R}\cup\set{+\infty}$ at position
$\bar{x}$ is a set valued mapping given by
\begin{equation}
  \label{eq:4}
  \partial \varphi\left(\bar{x}\right) \coloneqq
  \set{ x^{*} \middle|\
    \varphi \left( x \right) - \varphi \left( \bar{x} \right)
    \geqslant \langle x^* , x - \bar{x} \rangle,\ \forall x}\enspace{}.
\end{equation}
Its elements are called subgradients. Without further requirements on $\varphi$
this set may contain infinitely many elements or be empty. A common example is
the subgradient of the absolute value function where a simple computation shows
that
\begin{equation}
  \label{eq:5}
  \partial \left( \left|\:\cdot\:\right| \right)\left( x \right) =
  \begin{cases}
        \set{-1},          &x<0\enspace{},\\
        \left[-1,1\right], &x=0\enspace{},\\
        \set{1},           &x>0\enspace{}.
  \end{cases}
\end{equation}
For strictly concave functions the subdifferential is always empty. On the other
hand, convex functions always have a least one subgradient in the interior of
their domain. Subdifferentials exhibit many properties of usual derivatives. One
of the most important properties of subdifferentials is for example that
$0\in\partial \varphi\left(\bar{x}\right)$ is a necessary condition for
$\bar{x}$ being a minimiser of $\varphi$.\par
Robust regularisers involving the $\ell_{1}$ norm are quite common in
variational image analysis models. Their optimisation often leads to subproblems
of the kind
\begin{equation}
  \label{eq:1}
  \argmin_{x\in\mathbb{R}^{n}}\set{\norm{x}_{1} +
    \frac{\lambda}{2}\norm{x - b}_{2}^{2}}
\end{equation}
with a positive parameter $\lambda$ and an arbitrary vector $b$. A closed form
solution $x^{*}$ can be derived in terms of the well known soft shrinkage
operator: $x^{*}=\shrink (b, \frac{1}{\lambda})$, where
\begin{equation}
  \label{eq:300}
  \shrink \left( y, \alpha \right) \coloneqq
  \begin{cases}
    y - \alpha, &y > \alpha\enspace{},\\
    0,          &y\in\left[-\alpha,\alpha\right]\enspace{},\\
    y + \alpha, &y < -\alpha\enspace{}.
  \end{cases}
\end{equation}
For vector valued arguments, the shrinkage operator is applied componentwise.
Unfortunately, the $\ell_{1}$ norm is not rotationally invariant and promotes in
many applications undesired structures parallel to the coordinate axes. A
possible workaround consists in adapting the considered models such that we are
lead to tasks of the form
\begin{equation}
  \label{eq:2}
  \argmin_{x\in\mathbb{R}^{n}}\set{\norm{x}_2 + \frac{\lambda}{2}\norm{x - b}_2^2}
\end{equation}
with $\lambda>0$. Here, the closed form solution can be expressed in terms of
the generalised shrinkage operator.
\begin{definition}[Generalised Shrinkage]
  \label{thm:1}
  Let $b$ be a vector in $\mathbb{R}^n$ and $\lambda >
  0$, then we define the \emph{generalised shrinkage operator} as
  \begin{equation}
    \begin{split}
      \gshrink\left(b,\lambda\right)
      &\coloneqq
      \max\left(\norm{b}_2-\lambda,0\right)\frac{b}{\norm{b}_2}
      \phantom{:}=
      \begin{cases}
            b - \frac{\lambda}{\norm{b}}_2 b,
            &\text{if}\ \norm{b}_2>\lambda\enspace{},\\
            0, 
            &\text{else}\enspace{},
      \end{cases}
    \end{split}
    \end{equation}
  where we adopt the convention $0\cdot\frac{0}{0}=0$.
\end{definition}
The solutions $x^{*}$ of \eqref{eq:2} can then be expressed in term of this
generalised shrinkage. It holds
$x^{*}=\gshrink\left(b,\frac{1}{\lambda}\right)$. The proof is lengthy but not
difficult. One has to find those $x^{*}$ for which 0 is a subgradient of the
cost function. This can be done by discerning the cases $\lambda \norm{b}_2 > 1$
and $\lambda \norm{b}_2 \leqslant 1$.
%
\section{The Bregman Framework}
\label{sec:3}
%
We begin with recalling the standard Bregman iteration as developed by Osher
\emph{et al.}\ \cite{Osher2005}. Furthermore, we present an alternative but
equivalent formulation of the Bregman algorithm which has been discussed in
\cite{Goldstein2009}. We make use of this formulation as it simplifies the proof
of a convergence assertion and for describing the SBM introduced by Goldstein
and Osher \cite{Goldstein2009}. As indicated, the SBM will be the basis for our
OF algorithms. Let us emphasise that many approaches to the Bregman framework
exist in the literature. Bregman himself \cite{Bregman1967} wanted to describe
non-orthogonal projections onto convex sets. Further research in that direction
can for example be found in \cite{BB1997}. In \cite{E2009,Setzer2009} a certain
number of equivalences between different optimisation techniques are discussed.
They allow us to interpret the Bregman algorithms by means of conjugate duality
for convex optimisation. Thus, the Bregman framework may also be interpreted as
a splitting scheme or a primal dual algorithm. The presentation in this work
relies more on similarities between constrained and unconstrained convex
optimisation problems. The convergence theory that we will recall and present
here is based on results in
\cite{Brune2009,Burger2006,BRH07,Cai2009a,Goldstein2009,Osher2005}. The authors
of these works employed different mathematical settings. Some of the results
require Hilbert spaces, others are stated in rather general vector spaces only
equipped with semi-norms. We unify the results here within one framework using
finite dimensional, normed vector spaces. This set-up allows the usage of a
common set of requirements and to clarify relations between the different works.
While doing this, we also add some new results to the SBM.\par
We note that this mathematical setting suffices for the typical application in
computer vision where one ultimately needs to resort to a discretised
problem.\par
Let us now introduce the mathematical formalism that we will require in the
forthcoming sections. One of the central concepts behind the Bregman iteration
is the Bregman divergence. It has been presented by Bregman in 1967
\cite{Bregman1967}, where it has been used to solve convex optimisation problems
through non-orthogonal projections onto convex sets.\par
\begin{definition}[Bregman Divergence]
      \label{thm:2}
      The \emph{Bregman divergence}
      $D_\varphi^p \colon \mathbb{R}^{n} \times \mathbb{R}^{n} \to \mathbb{R}$
      of a proper convex function $\varphi$ is defined as $D_\varphi^p\left(
            x,y \right)\coloneqq \varphi \left( x \right) - \varphi
      \left( y \right) - \langle p, x - y \rangle$. Thereby
      $p$ is a subgradient of $\varphi$ at $y$.
\end{definition}
The aim of the Bregman iteration is to have a formulation that can handle convex
non-differentiable cost functions and that avoids ill conditioned formulations.
To illustrate the main idea, one may consider e.g.\ the following optimisation
problem:
\begin{equation}
  \label{eq:101}
  x^{(k+1)} =
  \argmin_{x\in \mathbb{R}^{n}}\ \set{D_\varphi^{p}\left( x , x^{(k)} \right)
    + \iota_{\set{0}}\left( Ax-b \right)}
\end{equation}
where $\iota_{S}$ is the indicator function of the set $S$, i.e.\ $\iota_{S}(x)$
is 0 if $x\in S$ and $+\infty$ else. In case the linear system $Ax=b$ has
multiple solutions or if it has a very large system matrix, then it might be
difficult to determine the iterates $x^{(k)}$. Therefore, one may reformulate
\eqref{eq:101} in terms of a regularised and unconstrained problem
\begin{equation}
  \label{eq:102}
  x^{(k+1)} =
  \argmin_{x\in \mathbb{R}^{n}}\set{D_\varphi^{p} \left( x , x^{(k)} \right) +
    \lambda \norm{Ax-b}_{2}^{2}}
\end{equation}
with some fixed $\lambda>0$ to approximate \eqref{eq:101}. This iterative
strategy motivates Definition~\ref{thm:3}, which coincides with the formulation
found in \cite{Goldstein2009,Osher2005,Yin2007}.\par
Let us note that in \cite{Osher2005,Yin2007} the Bregman iteration has been
formulated as a method for minimising convex functionals of the form
$J\left(u\right) + H\left(u\right)$. However, the convergence theory presented
below that is derived from
\cite{Brune2009,Burger2006,BRH07,Cai2009a,Goldstein2009,Osher2005} states that
the iterates converge towards the solution of a constrained formulation.
Therefore, we define the algorithm from the beginning on as a method for solving
constrained optimisation problems.
In the following we silently assume that $J$ and $H$ are always two proper
convex functions defined on the whole $\mathbb{R}^n$. Further, $H$ will be a
non-negative differentiable function with $0 = \min_{u}\set{H\left(u\right)}$
and this minimum is reached at some point in $\mathbb{R}^{n}$.\par
\begin{definition}[Bregman iteration]
  \label{thm:3}
  The \emph{Bregman iteration} of the constrained optimisation problem
  \begin{equation}
    \label{bh-constrainedproblem}
    \argmin_{u\in\mathbb{R}^{{n}}} \set{J\left(u\right) +
      \iota_{\set{0}}\left( H\left( u \right) \right)}
  \end{equation}
  is given by:
  \begin{enumerate}
    \item Choose $u^{(0)}$ arbitrarily, $\lambda >0$ and
          $p^{(0)} \in \partial J\left(u^{(0)}\right)$.
    \item Compute iteratively
          \begin{equation}
            \label{bh-bregmaniter}
            u^{(k+1)} = \argmin_{u\in\mathbb{R}^{{n}}}
            \set{D_{J}^{p^{(k)}} \left( u,u^{(k)} \right) +
              \lambda H \left( u \right)}\\
          \end{equation}
          where we have $p^{(k)} \in \partial J \left( u^{(k)} \right)$
          until a fixed-point is reached.
  \end{enumerate}
\end{definition}
From our assumptions on $J$ it follows that it has at least one subgradient at
every point. Thus $p^{(k)}$ always exists, but it is not necessarily unique. In
general settings the computation of a subgradient may not always be simple. The
following result from \cite{Yin2007} provides a comfortable strategy to obtain a
single specific subgradient.
\begin{proposition}
  \label{bh-BregmanIterDef2}
  If $H$ is differentiable, the second step of the Bregman iteration from
  Definition \ref{thm:3} becomes
  \begin{equation}
    \label{bh-bregmaniter1}
    \begin{split}
          u^{(k+1)} &=
          \argmin_{u\in\mathbb{R}^{n}} \set{D_J^{p^{(k)}}\left(u,u^{(k)}\right)
            + \lambda H \left( u \right)}
          \enspace{},\\
          p^{(k+1)} &= p^{(k)} - \lambda \nabla H\left(u^{(k+1)}\right)
          \enspace{}.
    \end{split}
  \end{equation}
  \begin{proof}
        It suffices to show that $p^{(k)} - \lambda \nabla H(u^{(k+1)})$ is a
        subgradient of $J$ at position $u^{(k+1)}$. The definition of the
        iterates $(u^{(k)})_{k}$ implies that $u^{(k+1)}$ is a minimiser of
        $D_J^{p^{(k)}} (u,u^{(k)} ) + \lambda H \left( u \right)$.
        Expanding the definition of the Bregman divergence and removing the
        constant terms, we see that $0$ is a subgradient of
        $J(u) - \langle p^{(k)}, u \rangle + \lambda H ( u )$ at position
        $u^{(k+1)}$. Since the subdifferential of a sum coincides with the sum
        of the subdifferentials it follows that there must exist
        $p^{(k+1)} \in \partial J (u^{(k+1)})$ that fulfils the equation
        \begin{equation}
              \label{eq:6}
              0 = p^{(k+1)} - p^{(k)} + \lambda \nabla H\left(u^{(k+1)}\right)
              \enspace{}.
        \end{equation}
  \end{proof}
\end{proposition}
Although one could basically use any subgradient of $J$ at $u^{(k+1)}$, the
previous proposition gives us a convenient way of finding a specific one that is
easy to obtain. This makes the computation of the iterates much simpler and
improves the overall speed of the algorithm.\par
Our next goal is to analyse the convergence behaviour of the Bregman iteration
given in Definition~\ref{thm:3}, and to show that its iterates obtained from
\begin{equation}
  \label{eq:16}
  \argmin_{u\in\mathbb{R}^{n}}\ \set{D_J^{p^{(k)}}\left(u,u^{(k)}\right) + \lambda H\left(u\right)}
\end{equation}
converge towards a solution of
\begin{equation}
  \label{bh-consideredproblem}
  \argmin_{u\in\mathbb{R}^{n}} \set{J\left(u\right) +
    \iota_{\set{0}} \left( H\left(u\right) \right)}\enspace{}.
\end{equation}
For well-posedness reasons we will assume that \eqref{bh-consideredproblem} as
well as the iterative formulation in \eqref{eq:16} are always solvable. If these
requirements are not feasible, then our iterative strategy cannot be carried out
or fails to converge. We emphasise that the existence of a solution of either
\eqref{eq:16} or \eqref{bh-consideredproblem} cannot always be deduced from the
existence of the other formulation. Even if $H(u)=0$ cannot be fulfilled, it
might still be possible that all iterates in \eqref{eq:16} exist.\par
The results compiled by the Propositions~\ref{bh-bregmantheoconv},
\ref{bh-bregmantheoconvfoo}, and~\ref{bh-convsolH}, as well as
Corollary~\ref{bh-rembregdist} were already discussed in \cite{Osher2005}.
There, the authors discussed iterative regularisation strategies in the space of
functions with bounded variation. Furthermore, they established the link between
their algorithm and the Bregman divergence. The proofs from \cite{Osher2005} for
the variational setting carry over verbatim to the finite dimensional set-up
that we use within this work. Thus, we just recall the statements without
proofs.
\begin{proposition}
  \label{bh-bregmantheoconv}
  The sequence $(H(u^{(k)}))_k$ is monotonically decreasing. We have
  for all~$k$:
  \begin{equation}
    H\left(u^{(k+1)}\right) \leqslant H\left(u^{(k)}\right)
  \end{equation}
  and strict inequality when $D_{J}^{p^{(k)}}(u^{(k+1)},u^{(k)})$ is positive.
\end{proposition}
In this context we remark that the Bregman divergence is always non-negative for
convex $J$ and if $J$ is even strictly convex, then $D_{J}^{p}(x, y) = 0$ can
only hold if and only if $x=y$.
\begin{proposition}
  \label{bh-bregmantheoconvfoo}
  We have for all $\lambda > 0$
  \begin{multline}
    \label{bh-convbregman-2} D^{p^{(k)}}_J\left(u,u^{(k)}\right) +
    D^{p^{(k-1)}}_J\left(u^{(k)},u^{(k-1)}\right) \\ -
    D^{p^{(k-1)}}_J\left(u,u^{(k-1)}\right) \leqslant
    \lambda \left( H\left(u\right) - H\left(u^{(k)}\right) \right).
  \end{multline}
\end{proposition}
\begin{corollary}
  \label{bh-rembregdist}
  For the particular choice $u = \tilde{u}$, where $\tilde{u}$ is a solution of
  $H\left(u\right)=0$ we immediately get:
  \begin{equation}
    0 \leqslant D^{p^{(k)}}_J\left(\tilde{u},u^{(k)}\right)
    \leqslant D^{p^{(k-1)}}_J\left(\tilde{u},u^{(k-1)}\right)\enspace{}.
  \end{equation}
\end{corollary}
One can easily infer from the above assertions, that for strictly convex $J$,
the iterates converge towards a solution of $H(u)=0$. The next proposition
gives an estimate how fast this convergence is, and it shows that the strict
convexity is in fact not necessary. Let us note that the proof of it relies on
the Propositions~\ref{bh-bregmantheoconv} and \ref{bh-bregmantheoconvfoo}.
\begin{proposition}\label{bh-convsolH}
  If $\tilde{u}$ is a solution of $H(u)=0$ and if
  $D_J^{p^{(0)}}\left(\tilde{u},u^{(0)}\right) < +\infty$ for some starting value
  $u^{(0)}$ then one obtains for all $\lambda > 0$
  \begin{equation}
    \label{bh-convbregman-4}
    0 = H\left(\tilde{u}\right) \leqslant H\left(u^{(k)}\right)
    \leqslant
    \frac{D_J^{p^{(0)}}\left(\tilde{u},u^{(0)}\right)}{\lambda k}\enspace{}.
  \end{equation}
  Therefore, the iterates $u^{(k)}$ always converge towards a solution of $H$.
\end{proposition}
So far we have seen that the iterates converge towards a solution of $H$. But at
this point we do not know whether this solution also minimises our cost function
$J$. If $H$ has a unique solution, then the above theory is already
sufficient.\par
The following proposition states, that under certain assumptions, the iterates
that solve $H$ also minimise our cost function, even if $H$ has multiple
solutions. This highly important result was first pointed out in \cite{Yin2007},
where the authors analysed the convergence behaviour of the Bregman iteration
within the context of the basis pursuit problem. There, the authors analysed the
Bregman framework in a finite dimensional setting and further showed an
interesting relationship to augmented Lagrangian methods.
\begin{proposition}
  \label{bh-bregmansolvesconst}
  Assume there exists $u^{(0)}$ such that it is possible to choose $p^{(0)}=0$
  in \eqref{bh-bregmaniter1}. Furthermore, assume that
  $H\left(u\right) = h\left( Au-b \right)$ where $A$ is some matrix, $b$ an
  arbitrary vector and $h$ is a differentiable non-negative convex function that
  only vanishes at $0$. If an iterate $u^{(k)}$ fulfils $H(u^{(k)}) = 0$, i.e.\
  it solves $Au=b$, then that iterate is also a solution of the constrained
  optimisation problem of \eqref{bh-constrainedproblem}.
\end{proposition}
Concerning the proof of Proposition~\ref{bh-bregmansolvesconst} from
\cite{Yin2007}, let us note that this proposition requires that $h$ vanishes
only at $0$. However, the linear system $Au=b$ can have multiple solutions. Thus
$H$ can have multiple solutions, too. The requirement that $h$ only vanishes at
$0$ is essential in the proof, as it enforces that every zero of $H$ solves the
linear system.\par
We conclude that convergence is guaranteed if $H$ has the form described in
Proposition~\ref{bh-bregmansolvesconst} and if we can choose $u^{(0)}$ such that
$p^{(0)} = 0$ is a valid subgradient. The latter requirement is in fact rather
weak such that only the former is of importance. For the formulation of the OF
problems, these conditions will fit naturally into the modelling.\par
We will now focus on the special case of interest for us that
$h\left(x\right)\coloneqq\norm{x}_2^2$. In that case it is possible to derive an
estimate for the error at each iteration step. Such a result was presented in
\cite{BRH07}, where the authors discussed the convergence behaviour of the
Bregman iteration in the context of inverse scale methods for image restoration
purposes. Their setting included a variational formulation and used $L^{p}$
spaces as well as the space of functions of bounded variation. Furthermore, they
had to formulate certain convergence results in terms of the weak-* topology.
The usage of finite dimensional settings allows a more consistent formulation.
In our mathematical setting, the proof can be done analogously to the one in
\cite{BRH07}.\par
In order to prepare the presentation of the SBM formulation we consider in the
following a more concrete optimisation task. Therefore, assume now that
$A$ is a given $m\times n$ matrix and $b$ a known vector in
$\mathbb{R}^{n}$. The problem that we consider is
\begin{equation}
  \label{bh-convbregman-7}
  \argmin_{u\in\mathbb{R}^{n}}\set{J\left(u\right) +
    \iota_{\set{0}}\left( \frac{1}{2} \norm{Au-b}^2_2 \right)}\enspace{}.
\end{equation}
Proposition~\ref{bh-BregmanIterDef2} implies that we have the following
algorithm:
\begin{equation}
  \label{bh-convprob1}
  \begin{split}
    u^{(k+1)} &= \argmin_{u\in\mathbb{R}^{n}}\set{D_{J}^{p^{(k)}}(u,u^{(k)}) +
      \frac{\lambda}{2}\norm{Au-b}^{2}_{2}},\\
    p^{(k+1)} &= p^{(k)} + \lambda A^{T}\left( b - Au^{(k+1)} \right)\enspace{}.
  \end{split}
\end{equation}
For technical reasons we will continue to assume that it is possible
to choose $u^{(0)}$ such that $p^{(0)}=0$ can be used. This can always be done
as long as $J$ has a minimum at some finite point in our framework. It is
useful to consider the following two definitions, which stem from \cite{BRH07}.
\begin{definition}[Minimising Solution]
  \label{bh-def-minsol}
  A vector $\tilde{u}$ is called $J$ \emph{minimising solution} of $Au=b$, if
  $A\tilde{u}=b$ and $J\left(\tilde{u}\right)\leqslant J\left(v\right)$ for all
  other $v$ that fulfil $Av=b$.
\end{definition}
\begin{definition}[Source Condition]
  \label{bh-def-sourcond}
  Let $\tilde{u}$ be a $J$ minimising solution of $Au=b$. We say $\tilde{u}$
  satisfies the \emph{source condition} if there exists an $\omega$ such that
  $A^{T}\omega \in \partial J\left(\tilde{u}\right)$.
\end{definition}
The source condition can, in a certain sense, be interpreted as an additional
regularity condition that we impose on the solution. Not only do we require
that the minimising solution has a subgradient, we even want that there exists
a subgradient that lies in the range of $A^T$. Requirements like this
are a frequent tool in the analysis of inverse problems.
The next theorem adopted from \cite{BRH07} shows that it is possible to
give an estimate for the error if this source condition holds.
\begin{theorem}
  \label{bh-theorem-convspeedBreg}
  Let $\tilde{u}$ be a solution of $Au=b$ minimising $J$, and assume that the
  source condition holds, i.e. there exists a vector
  $\xi \in \partial J\left(\tilde{u}\right)$ such that $\xi = A^{T}q$ for some
  vector $q$. Furthermore, assume that it is possible to choose $u^{(0)}$
  such that $p^{(0)}=0$ is a subgradient of $J$ at $u^{(0)}$.
  Then we have the following
  estimation for the iterates $u^{(k)}$ of \eqref{bh-convprob1}:
  \begin{equation}\label{bh-convbregman-8}
    D^{p^{(k)}}_J \left(\tilde{u},u^{(k)}\right)
    \leqslant
    \dfrac{\norm{q}^2_2}{2\lambda k}\quad \forall k \in \mathbb{N}^{*}\enspace{}.
  \end{equation}
\end{theorem}
The result of the following proposition can be found in \cite{Yin2007}.
\begin{proposition}[Alternative formulation]
  \label{bh-propsimpl}
  The normal Bregman iteration for solving the constrained optimisation problem
  \begin{equation}\label{bh-simple-constprob}
    \argmin_{u\in\mathbb{R}^{n}}\set{J\left(u\right) +
      \iota_{\set{0}}\left( \frac{1}{2} \norm{Au-b}^{2}_{2} \right)}
  \end{equation}
  can also be expressed in the following iterative form:
  \begin{equation}
    \label{bh-bregmanitersimp}
    \begin{split}
      u^{(k+1)} &\coloneqq
      \argmin_{u\in\mathbb{R}^{n}}\set{J\left(u\right) +
        \frac{\lambda}{2}\norm{Au-b^{(k)}}^{2}_{2}}\enspace{},\\
      b^{(k+1)} &\coloneqq b^{(k)} + b - Au^{(k+1)}\enspace{}.
    \end{split}
  \end{equation}
  if we set $b^{(0)} \coloneqq b$ and choose $u^{(0)}$ such that $p^{(0)}=0$.
\end{proposition}
Because of the equivalence of the two formulations the iterates
given by this alternative Bregman algorithm have the same properties as the ones
of the standard Bregman iteration. Thus, all the convergence results for the
standard set-up also apply in this case.\par
%
\subsection{The Split Bregman Method}
\label{sec:31}
%
The split Bregman method (SBM) proposed in \cite{Goldstein2009} extends the
Bregman iteration presented so far. It aims at minimising unconstrained convex
energy functionals.\par
While we mainly follow \cite{Goldstein2009} for the description of the
algorithm, we will also give some new results. We will for example discuss how
the convergence estimate of Brune \emph{et al.}\ \cite{Brune2009} can be applied
to the SBM.\par
The split Bregman formulation is especially useful for solving the following two
problems:
\begin{equation}
  \label{bh-sbm-1}
  \argmin_{u\in\mathbb{R}^{n}}\
  \set{\norm{\Phi\left(u\right)}_k + G\left(u\right)} \enspace{},\quad{}k=1,2
  \enspace{}.
\end{equation}
The function $\Phi$ is an affine mapping, i.e.\ $\Phi\left(u\right)=\Lambda u+b$
for some matrix $\Lambda$ and some vector $b$. $G$ should be a convex function
from $\mathbb{R}^n$ to $\mathbb{R}$. The difficulty in minimising these cost
functions stems from the fact that neither $\norm{\cdot}_1$, nor
$\norm{\cdot}_2$ are not differentiable in~0.\par
The basic idea behind the SBM is to introduce an additional variable that
enables us to separate the non-differentiable terms from the differentiable
ones. This is done by rewriting \eqref{bh-sbm-1} as a constrained optimisation:
\begin{equation}\label{bh-sbm-2}
  \argmin_{d,u\in\mathbb{R}^{n}}\set{\norm{d}_k + G\left(u\right) +
    \iota_{\set{0}}(d-\Phi\left(u\right))} \enspace{}.
\end{equation}
The previous section has shown us how to handle constrained optimisation tasks
of this kind. The main idea of SBM is to apply standard Bregman to
\eqref{bh-sbm-2}. In order to simplify the presentation, we employ the following
aliases:
\begin{align}\label{bh-sbm-3}
  \eta &\coloneqq \left(u,d\right)^T\; ,\\
  J\left(\eta\right)&\coloneqq\norm{d}_k + G\left(u\right)\; ,\\
  A\left( \eta \right)&\coloneqq d - \Lambda u\enspace{}.
\end{align}
Obviously $J$ is again a convex function and $A$ is a linear mapping.
Using the new notations, \eqref{bh-sbm-2} can be rewritten as
\begin{equation}\label{bh-sbm-4}
  \argmin_{\eta\in\mathbb{R}^{n}}\set{J\left(\eta\right) +
    \iota_{\set{0}}
    \left( \frac{1}{2}\norm{ A \left( \eta \right) - b }_2^2 \right)}\enspace{}.
\end{equation}
We assume at this point that it is possible to
choose $\eta^{(0)}$ such that $0$ is a subgradient of $J$ at $\eta^{(0)}$.
This is always possible if $G$ attains its minimum.\par
By applying the Bregman algorithm from Proposition~\ref{bh-propsimpl} one
obtains the following iterative procedure:
\begin{equation}
  \begin{split}
    b^{(0)} &= b\enspace{}, \\
    \eta^{(k+1)} &=
    \argmin_{\eta\in\mathbb{R}^{n}}\set{J\left(\eta\right) +
      \frac{\mu}{2} \norm{A\left(\eta\right) - b^{(k)}}_2^2},\\
    b^{(k+1)} &= b^{(k)} + b - A\left(\eta^{(k+1)}\right)
  \end{split}
\end{equation}
with $\mu>0$ being a constant positive parameter. Reintroducing the definitions
of $J$ and $A$ leads to a simultaneous minimisation of $u$ and $d$. Since such
an optimisation is difficult to perform, we opt for an iterative alternating
optimisation with respect to $u$ and $d$. Goldstein \emph{et al.}\ \cite{Goldstein2009}
suggest to do a single sweep. In this paper we allow a more flexible handling
and allow up to $M$ alternating optimisations. All in all, we have to solve for
$j=1,\ldots,M$:
\begin{align}
  \argmin_{u\in\mathbb{R}^{n}}
  &\set{ G\left(u\right) + \frac{\mu}{2} \norm{d^{(k,j)} - \Lambda u - b^{(k)}}_2^2}
    \enspace{},\\
  \argmin_{d\in\mathbb{R}^{n}}
  &\set{\norm{d}_1 + \frac{\mu}{2} \norm{d - \Lambda u^{(k,j+1)} - b^{(k)}}_2^2}
    \enspace{}.
\end{align}
The first optimisation step depends largely on the exact nature of $G$. As a
consequence one cannot make any general claims about it. We just note that for
the case where $G\left(u\right) = \norm{Cu-f}_2^2$ for some matrix $C$ and a
vector $f$, the cost function becomes differentiable and the minimiser can be
obtained by solving a linear system of equations with a positive semi-definite
matrix. If either $C$ or $\Lambda$ has full rank, then the system matrix will
even be positive definite. This will especially be true for the upcoming
applications to optic flow. The second optimisation has a closed form solution
in terms of shrinkage operations. The solution is given by
\begin{equation}
  d^{(k,j+1)} =
  \shrink\left( \left(\Lambda u^{(k,j+1)} + b^{(k)}\right) , \frac{1}{\mu}\right)
\end{equation}
where the computation is done componentwise. If we replace the $\ell_{1}$ norm
by the Euclidean norm, then we have to resort to the generalised shrinkage
operator and the solution is given by
\begin{equation}
  d^{(k,j+1)} = \gshrink\left( \Lambda u^{(k,j+1)} + b^{(k)} , \frac{1}{\mu} \right)
  \enspace{}.
\end{equation}
The detailed formulation of the SBM with $N$ iterations and $M$ alternating
minimisation steps for solving \eqref{bh-sbm-1} is depicted in
Algorithm~\ref{bh-algor1}.\par
\begin{algorithm*}
  \dontprintsemicolon
  \KwData{$\Lambda$, $b$, $G$, $N$, $M$}
  \KwResult{$u^{(N)}$ and $d^{(N)}$ minimising \eqref{bh-sbm-1}}
  \KwInit{$u^{(0)}$ such that $0 \in \partial G\left(u^{(0)}\right)$, 
    $d^{(0)}=0$, $b^{(0)}=b$}
  \For{$k=0$ \KwTo $N-1$}{
    $u^{(k,0)} = u^{(k)}$ and $d^{(k,0)} = d^{(k)}$\;
    \For{$j=0$ \KwTo $M-1$}{
      $\displaystyle u^{(k,j+1)} = \argmin_{u}\ \set{G\left(u\right) + 
        \frac{\mu}{2} \norm{d^{(k,j)} - \Lambda u - b^{(k)}}_2^2}$\;
      $\displaystyle d^{(k,j+1)} = \argmin_{d}\ \set{\norm{d}_1 + 
        \frac{\mu}{2} \norm{d - \Lambda u^{(k,j+1)} - b^{(k)}}_2^2}$\;
    }
    $u^{(k+1)} = u^{(k,M)}$ and $d^{(k+1)} = d^{(k,M)}$\;
    $b^{(k+1)} = b^{(k)} + b - d^{(k+1)} + \Lambda u^{(k+1)}$}
  \caption[Split Bregman algorithm with alternative Bregman iteration.]{Split
    Bregman algorithm for $N$ iterations with $M$ alternating minimisation steps
    based upon the alternative form of the Bregman iteration.}
  \label{bh-algor1}
\end{algorithm*}
Since the SBM relies on the Bregman iteration, it is clear that all the related
convergence results also hold for the SBM. Especially
Theorem~\ref{bh-theorem-convspeedBreg} gives us an estimate for the convergence
speed if certain regularity conditions are met. In the following, we would like
to analyse if these conditions can be fulfilled for the SBM. We are going to
consider the following problem
\begin{equation}
  \label{bh-spBeq1}
  \argmin_{u} \set{\sum_{k=1}^N \norm{A_{k} u + b_k }_2 + \frac{\lambda}{2} \norm{B u -c}_2^2}
\end{equation}
where $A_k : \mathbb{R}^n\to\mathbb{R}^m$, $B : \mathbb{R}^n\to\mathbb{R}^l$ are matrices,
$b_k\in\mathbb{R}^m$, $c\in\mathbb{R}^l$ some vectors and $\lambda > 0$ a real-valued
parameter. This model represents a generic formulation that also includes all
forthcoming OF models. Thus, all statements concerning this model are
automatically valid for our OF methods, too. The corresponding split Bregman
algorithm of \eqref{bh-spBeq1} solves
\begin{equation}
  \label{bh-spBeq2}
  \begin{split}
    &\argmin_{u,d_1,\ldots,d_N} \set{\sum_{k=1}^N \norm{d_k}_2 + \frac{\lambda}{2} \norm{Bu-c}_2^2}\\
    &\quad\text{such that}\quad\quad\sum_{k=1}^N\norm{d_k - A_k u - b_k }_2^2 = 0\enspace{}.
  \end{split}
\end{equation}
Note that the necessary conditions for the application of the split Bregman algorithm
are met. The cost function attains its global minimum for $d_k = 0$ for all $k$ and
$u$ a solution of $B^TBu=B^Tc$. The constraining condition obviously also has a solution.
Let us now define the matrix $\Lambda : \mathbb{R}^{n+Nm}\to\mathbb{R}^{Nm}$ by
\begin{equation}
  \Lambda\coloneqq
  \begin{bmatrix}
    -A_1   & I & 0 & \ldots & 0 \\
    -A_2   & 0 & I & \ldots & 0 \\
    \vdots & \vdots & \vdots & \ddots & \vdots \\
    -A_N   & 0 & 0 & \ldots & I \\
  \end{bmatrix}
\end{equation}
where $I$ is the identity matrix in $\mathbb{R}^m$. If we define further
$\tilde{b}=\left(b_1,b_2,\ldots,b_N\right)^T$ then \eqref{bh-spBeq2} can be
rewritten as
\begin{equation}
  \label{bh-spBeq3}
  \begin{split}
    &\argmin_{u,d_{1},\ldots,d_{N}}
    \biggl\{
    \underbrace{\sum_{k=1}^N \norm{d_k}_2 + 
      \frac{\lambda}{2} \norm{Bu-c}_2^2}_{\eqqcolon 
      J\left(u,d_1,\ldots,d_N\right)}\biggr\}\\
    &\quad\quad\text{such that}\quad\quad
    \norm{\Lambda\left(u,d_1,\ldots,d_N\right)^T - \tilde{b} }_2^2 = 0\enspace{}.
  \end{split}
\end{equation}
Now assume that we have found $u = (\tilde{u}, \tilde{d_1}, \ldots,
\tilde{d_N} )^T$, a $J$ minimising solution of $\Lambda x = \tilde{b}$. In
order to apply Theorem~\ref{bh-theorem-convspeedBreg} we need to know how
$\partial J\left( \tilde{u} , \tilde{d_1} , \ldots , \tilde{d_N} \right)$ looks
like. So assume $\left(w, w_1, \ldots, w_N\right)^T$ is a subgradient. By
definition we must have for all $d_k,\ k=1,\ldots,N$ and all $u$:
\begin{multline}
  \label{eq:7}
  \sum_{k=1}^N \left( \norm{d_k}_2  - \norm{\tilde{d_k}}_2 \right) + 
  \frac{\lambda}{2}(\norm{Bu-c}_2^2 - \norm{B\tilde{u}-c}_2^2) \\
  \geqslant \langle w , u - \tilde{u} \rangle + \sum_{k=1}^N \langle w_k , d_k - \tilde{d_k} \rangle\enspace{}.
\end{multline}
Since this must hold for all possible choices, it must hold especially for $d_k
= \tilde{d_k}$ with $k=1,\ldots , N$. But then we see that $w$ must be a
subgradient of $\frac{\lambda}{2}\norm{Bu-c}_2^2$ at $\tilde{u}$. Setting
$u=\tilde{u}$ and all but one $d_k$ to $\tilde{d_k}$ yields in the same way that
every $w_k$ must be a subgradient of $\norm{d}_2$ at $\tilde{d_k}$. It follows
that we have the following representation
\begin{align}
  w &= \lambda B^T\left(B\tilde{u}-c\right)\enspace{}, \label{eq:8}\\
  w_k &= \frac{\tilde{d_k}}{\norm{\tilde{d_k}}_2}\ \text{for}\ k=1,\ldots,N\enspace{}. \label{eq:9}
\end{align}
We assume here that all $\tilde{d_k}$ are different from $0$. If this is not the
case, then the choice of the subgradient is not unique anymore and would
complicate the following discussion. Theorem~\ref{bh-theorem-convspeedBreg}
requires that there is a vector
$\gamma\coloneqq\left(\gamma_1,\ldots,\gamma_N\right)^T$ such that
$\Lambda^T\gamma \in\partial J\left(\tilde{u}, \tilde{d_1}, \ldots,
\tilde{d_N}\right)$. From the structure of the matrix $\Lambda^T$ we deduce that
the following conditions must be fulfilled
\begin{equation}
  \label{bh-spBregeq4}
  \sum_{k=1}^N A_k^T \frac{\tilde{d_k}}{\norm{\tilde{d_k}}_2}
  = \lambda B^T\left(c-B\tilde{u}\right)\enspace{}.
\end{equation}
If this relation holds for the minimising solution, then the estimate
given in Theorem~\ref{bh-theorem-convspeedBreg} also holds for the split Bregman
algorithm.\par
Let us close this section with two small remarks concerning the previous
results.\par
Should any of the $\tilde{d_k}$ be 0, then any vector with norm less or equal
than $1$ would be a valid subgradient of $\norm{d}_2$ at $\tilde{d_k}$. In that
case we gain additional degrees of freedom in the above formula which increases
the chances that it can be fulfilled.\par
The SBM still converges even if \eqref{bh-spBregeq4} is not fulfilled.
Theorem~\ref{bh-theorem-convspeedBreg} only gives an estimate for the
convergence speed, not for the convergence itself. The convergence is guaranteed
by Propositions~\ref{bh-bregmantheoconv},~\ref{bh-convsolH} and
Proposition~\ref{bh-bregmansolvesconst}. We refer to \cite{NF2014} for an
additional discussion on the necessity criteria to assert convergence. Further
convergence investigations under duality considerations are also exhibited in
\cite{YMO2013}. A discussion on convergence rates under strong convexity and
smoothness assumptions can also be found in \cite{GB2016,G2016}. These works
also include findings on optimal parameter choices.
%
\section{Optic Flow: The Setup}
\label{sec:4}
%
The purpose of this section is to present the OF models that are addressed in
this work. First we briefly consider basic model components. Then we summarise
the models that are of interest here, in a variational set-up as well as in the
discrete setting.\par
%
\subsection{Optical Flow Models}
\label{sec:41}
%
Let us denote by $\partial^{(k)} f$ the set of all partial derivatives of order
less or equal than $k$ of a given image from a sequence $f:\Omega \times T \to
\mathbb{R}$, where $\Omega$ is a subset of $\mathbb{R}^2$ representing the (rectangular)
image domain and $T\subseteq \mathbb{R}$ is a time interval. We restrict our
attention to grey value images. Extensions to colour images are possible, they
just render the proceeding more cumbersome and offer little insight into the
underlying mathematics. The aim of the OF problem is to determine the flow field
of $f$ between two consecutive frames at the moments $t$ and $t+1$. The two
components of this displacement field are $u, v\colon\Omega\to\mathbb{R}$.\par
The general form of a variational model that determines the unknown displacement
field $(u,v)$ as the minimiser of an energy functional can then be written as
\begin{equation}
  \label{bh-energyfunctional}
  \argmin_{(u,v)} \set{\int_\Omega D\left(u,v\right) +
    \lambda S\left(\nabla u,\nabla v\right) \mathrm{d}y \mathrm{d}y}\enspace{}.
\end{equation}
Thereby, $D(u,v)$ denotes a data confidence term (or just data term),
while the so-called smoothness term $S\left(\nabla u, \nabla v \right)$
regularises the energy, and where $\lambda > 0$ is a regularisation parameter.
The operator $\nabla$ corresponds as usual to the gradient. Such variational
formulations have the advantages that they allow a transparent modelling, and
that the resulting flow fields are dense.\par
We employ a modern approach that combines the two following model assumptions:\par
\emph{Grey value constancy:} One assumes here that the following equality holds
\begin{equation}
  \label{bh-GreyValConst}
  f\left(x+u\left(x,y\right),y+v\left(x,y\right),t+1\right) = f\left(x,y,t\right)\enspace{}.
\end{equation}
Surprisingly, this assumption is relatively often fulfilled when the
displacements remain small enough. Unfortunately we have two unknowns but only
one equation. Thus, there is no chance to recover the complete displacement
field based on this equation alone. This problem is known in the literature as
the aperture problem.\par
\emph{Gradient constancy:} Assuming that the spatial gradient of $f$ remains
constant leads to the following equation
\begin{equation}
  \label{bh-GradValConst}
  \nabla f\left(x+u ,y+v ,t+1\right) = \nabla f\left(x,y,t\right)\enspace{}.
\end{equation}
Here we have two unknowns and two equations. As a consequence the aperture
problem is not always present. This assumption is of interest as it remains
fulfilled when the image undergoes global illumination changes, whereas the grey
value constancy does not.\par
Our constancy assumptions represent nonlinear relationships between the data $f$
and the flow field $(u,v)$. As a remedy we assume that all displacements are
small. In this setting we may approximate the left hand side of the equations
above by their corresponding first order Taylor expansions. Then
\eqref{bh-GreyValConst} becomes
\begin{equation}
  \label{eq:11}
  f_x u + f_y v + f_t = 0
\end{equation}
and \eqref{bh-GradValConst} becomes
\begin{equation}
  \label{eq:12}
  \begin{split}
    f_{xx} u + f_{xy} v + f_{xt} &= 0\enspace{}, \\
    f_{xy} u + f_{yy} v + f_{yt} &= 0
  \end{split}
\end{equation}
where the indices designate the partial derivatives with respect to the
corresponding variables. Deviations of the lefthand side from 0 can be
considered as errors and will be penalised in our models. Making use of a weight
$\gamma \geq 0$, interesting combinations of these models can be found in
Table~\ref{bh-table-terms-1}.\par
\begin{table}[t]
  \centering
  \caption{Summary of data and smoothness terms.}
  \label{bh-table-terms-1}
  \begin{tabular*}{0.75\linewidth}{@{\extracolsep{\fill} }ll} \toprule
    Term & Definition \\ \cmidrule{1-1}\cmidrule{2-2}
    $D_1$
         &
           $\displaystyle{
           \left( f_x u + f_y v + f_t \right)^2
           +
           \gamma
           \left( \left(
           f_{xx} u + f_{xy} v + f_{xt} \right)^2 + \left( f_{xy} u + f_{yy} v + f_{yt}
           \right)^2 \right)
           }$
    \\
    $D_2$
         &
           $\displaystyle{
           \abs{ f_x u + f_y v + f_t }
           +
           \gamma
           \left(
           \abs{ f_{xx} u + f_{xy} v + f_{xt} } + \abs{f_{xy} u + f_{yy} v + f_{yt}}
           \right)
           }$
    \\ \addlinespace
    $S_1$ & $\norm{\nabla u}^2 + \norm{\nabla v}^2$
    \\
    $S_2$ & $\norm{\nabla u} + \norm{\nabla v}$
    \\
    $S_3$ & $\sqrt{\norm{\nabla u}^2 + \norm{\nabla v}^2}$
    \\ \addlinespace \bottomrule
  \end{tabular*}
\end{table}
Using the notation $\norm{\nabla u}_2 \coloneqq \sqrt{ (\partial_x u)^2 +
(\partial_y u)^2 }$, we also address three smoothness terms of interest in OF
models, see $S_k$, $k=1,2,3$, in Table~\ref{bh-table-terms-1}.\par
The data term $D_1$ is optimal from a theoretical point of view because it is
convex and smooth. However, it is not robust with respect to outliers in the
data. The data term $D_2$ is more robust since the penalisation is
sub-quadratic. Its disadvantage is that it is not differentiable. The most
interesting smoothness terms are $S_2$ and $S_3$. Both are not differentiable,
but they offer a sub-quadratic penalisation. Furthermore $S_3$ is even
rotationally invariant. While $S_1$ is convex and differentiable and thus offers
attractive theoretical properties, the quadratic penalisation may cause an
oversmoothing of discontinuities in the motion field.\par
\begin{table}[t]
  \centering
  \caption{
    Summary of possible energy functionals:
    Horn and Schunck \cite{HS81},
    Pock \emph{et al.} \cite{Pock2008,Wedel2008,Pock2007},
    $OSB$: Optimal for split Bregman,
    $B$: Brox \emph{et al.} \cite{BBPW04}.\label{bh-table-energies}
  }
  \begin{tabular*}{0.75\linewidth}{@{\extracolsep{\fill} }lll}
    \toprule
    Energy & Data term & Smoothness term \\
    \cmidrule{1-1}\cmidrule{2-2}\cmidrule{3-3}
    \addlinespace
    Horn and Schunck & $D_1$ $(\gamma =0)$ & $S_1$ \\
    Pock \emph{et al.} & $D_2$ & $S_2$ \\
    \addlinespace
    $OSB$ (see Section~\ref{bh-sec-osb}) & $D_1$ & $S_3$ \\
    $B$ (see Section~\ref{bh-sec-broxmodel}) & $D_2$ & $S_3$ \\
    \addlinespace \bottomrule
  \end{tabular*}
\end{table}
Now that we have presented the smoothness and data terms, we can combine them to
different energy functionals. In Table \ref{bh-table-energies} we summarise the
possible choices and cite some references where these models have been
successfully applied.\par
%
\subsection{Algorithmic Aspects}
%
The following details are pre- and postprocessing steps that improve the
quality of our results. Most of these strategies are generic and many of them
are applied in various successful OF algorithms. We emphasise that they do not
infer with the Bregman framework that we use for the minimisation.\par
As usual for countless imaging applications we convolve each frame of our image
sequence $f$ with a Gaussian kernel with a small standard deviation in order to
deal with noise. For image sequences with large displacements, we follow
\cite{BBPW04} and embed the minimisation of our energy into a coarse-to-fine
multiscale warping approach. In all our experiments we set the scaling factor to
0.9. During warping, we employ a procedure from \cite{Wedel2008} where the
authors proposed to apply a median filter on the components
$\left(u^{c},v^{c}\right)$ that were obtained from the coarser grid. We point to
\cite{HB2011} for an analysis on the benefits of this strategy. Furthermore, we
disable the data term at occlusions. This can be achieved by multiplying the
data term with an occlusion indicator function $o \in \{0,1\}$, where $o(x,y)=0$
if a pixel is occluded and $o(x,y)=1$ if a pixel is visible. For the detection
of occlusions we follow the popular cross-checking technique from
\cite{cm92,pgpo94}. The occlusion handling is especially important for
approaches with a quadratic data term.\par
%
\section{Optical Flow: The Bregman Framework}
\label{sec:5}
%
In this section we elaborate on the formulation of the SBM for the considered OF
models. From an algorithmic point of view the most important questions have
already been answered. It remains to show that the optic flow models can be cast
into a form which is suitable for the application of the Bregman algorithms.
%
\subsection{The OSB model}
\label{bh-sec-osb}
%
First, we consider the model we denoted as OSB. A straightforward discretisation
yields
\begin{equation}
  \label{bh-osb-model}
  \argmin_{u,v} \set{\frac{\lambda}{2} D_1\left(u,v\right)
  + \sum_{i,j} \sqrt{\norm{\nabla u_{i,j}}_2^2 + \norm{\nabla v_{i,j}}_2^2}}
\end{equation}
where the summation goes over all pixel coordinates. Before we start applying
the Bregman algorithm let us have a look at the smoothness term first. In can be
reformulated in the following way
\begin{equation}
  \sum_{i,j}\ \sqrt{\norm{\nabla u_{i,j}}_2^2 + \norm{\nabla v_{i,j}}_2^2}
  =: \sum_{i,j} \norm{ \binom{\nabla u_{i,j}}{\nabla v_{i,j}} }_2\enspace{}.
\end{equation}
Thus, this model can also be written in the following more compact form
\begin{equation}
  \argmin_{u,v}\set{\frac{\lambda}{2} D_1\left(u,v\right) +
    \sum_{i,j} \norm{ \binom{\nabla u_{i,j}}{\nabla v_{i,j}} }_2}\enspace{}.
\end{equation}
The constrained formulation is now easily deduced. The best way to cast this
model into the SBM framework is to introduce slack variables $d^{u}_{i,j}$ and
$d^{v}_{i,j}$ for the non-differentiable smoothness term and to add and equality
constraint between the new variables and our flow field.
\begin{equation}
  \begin{split}
    &\argmin_{u,v}\ \set{\frac{\lambda}{2} D_1\left(u,v\right) + \sum_{i,j}
      \norm{\binom{d^u_{i,j}}{d^v_{i,j}}}_2}\\
    &\quad\quad\text{such that}\
    \frac{1}{2}\sum_{i,j} \norm{\binom{d^u_{i,j}}{d^v_{i,j}} - \binom{\nabla
        u_{i,j}}{\nabla v_{i,j}} }_2^2 = 0\enspace{}.
  \end{split}
\end{equation}
A straightforward reordering and grouping of all the involved terms leads us to
the following expression
\begin{equation}
  \begin{split}
        &\argmin_{u,v}\ \set{\frac{\lambda}{2} D_1\left(u,v\right) + 
          \sum_{i,j} \norm{\binom{d^u_{i,j}}{d^v_{i,j}}}_2}\\
        &\quad\quad\text{such that}\
        \frac{1}{2} \norm{\binom{d^u}{d^v} - \binom{\nabla u}{\nabla v} }_2^2 = 0
  \end{split}
\end{equation}
which is well suited for applying the Bregman framework. For convenience we have
grouped all variables $d_{i,j}^{u}$ and $d_{i,j}^{v}$ into large vectors
$d^{u}$, respectively $d^{v}$, while the vectors $\nabla u$ and $\nabla v$
contain the corresponding derivative information. The constraining condition
admits a trivial solution and thus, it does not pose any problem. The cost
function obviously has a minimum, too. The variables $d^u$ and $d^v$ act
independently of $u$ and $v$. Simply setting them all to $0$ and determining the
minimising $u$ and $v$ of $D_1$, by solving a least squares problem, yields the
desired existence of a minimiser. This implies that the cost function attains
its minimum and that there exists a point where $0$ is a subgradient. Note that
$D_1$ can always be minimised since we operate in a finite dimensional space,
where such problems are always solvable. It follows that the split Bregman
algorithm is applicable. Following the notational convention from
Section~\ref{sec:31}, we set
\begin{equation}
  \label{eq:13}
  \begin{split}
    \eta &= \left( u_{i,j}, v_{i,j}, d^{u}_{i,j}, d^{v}_{i,j} \right)^{\top}\enspace{},\\
    J\left( \eta \right)
    &= \frac{\lambda}{2}D_{1}(u,v) + \sum_{i,j} \norm{\binom{d^u_{i,j}}{d^v_{i,j}}}_2\enspace{}, \\
    A\left( \eta \right) &= \binom{d^u}{d^v} - \binom{\nabla u}{\nabla v}\enspace{},\\
    b &= 0\enspace{}.
  \end{split}
\end{equation}
The application of the SBM algorithm is now straightforward. In the alternating
optimisation steps the minimisation with respect to $(u,v)$ requires solving a
linear system of equations with a symmetric and positive definite matrix. We
refer to Section~\ref{linsystsec} for a proof. As mentioned in \cite{Goldstein2009} it
is enough to solve this system with very little accuracy. A few Gauß-Seidel
iterations are already sufficient. In \cite{YO2012}, the authors also discuss
the robustness of the Bregman approach with respect to inaccurate iterates and
provide a mathematically sound explanation. The minimisation with respect to
$d^u$ and $d^v$ can be expressed through shrinkage operations and does not pose
any problem. A detailed listing of the complete algorithm is given in
Algorithm~\ref{bh-L2L1Isoalgor}.
\begin{algorithm}
  \dontprintsemicolon
  \KwData{$F$, $F_x$, $F_y$, $f_t$, $f_{xt}$, $f_{yt}$,
    $\lambda$, $\gamma$, $\mu$, $N$, $M$}
  \KwResult{$u^{(N)}$ and $v^{(N)}$ minimising \eqref{bh-osb-model}}
  \KwInit{$u^{(0)}=v^{(0)}=0$, $d^{u,(0)}_{k,l}=0$, $b^{u,(0)}_{k,l}=0$,%
    $d^{v,(0)}_{k,l}=0$, $b^{v,(0)}_{k,l}=0$ for all $k$, $l$}
  \For{$i=0$ \KwTo $N-1$}{
    $u^{(i,0)} = u^{(i)}$, $v^{(i,0)} = v^{(i)}$\;
    $d^{u,(i,0)}_{k,l} = d^{u,(i)}_{k,l}$ and $d^{v,(i,0)}_{k,l} = d^{v,(i)}_{k,l}$
    for all $k$, $l$\;
    \For{$j=0$ \KwTo $M-1$}{
      $\displaystyle \left(u^{(i,j+1)},v^{(i,j+1)}\right) =
      \argmin_{u,v} \left\{ \frac{\lambda}{2} D_1\left(u,v\right) + \right.$\;
      \hfill$\displaystyle\left.\frac{\mu}{2}\left( \norm{d^{u,(i,j)} - 
                    \nabla u - b^{u,(i)}}_2^2 + \norm{d^{v,(i,j)} -
                    \nabla v - b^{v,(i)}}_2^2\right)\right\}$\;
      $\displaystyle \binom{d^{u,(i,j+1)}_{k,l}}{d^{v,(i,j+1)}_{k,l}} =
      \gshrink\left( \binom{\nabla u^{(i,j+1)}_{k,l}}{\nabla v^{(i,j+1)}_{k,l}} +
            \binom{ b^{u,(i)}_{k,l}}{ b^{v,(i)}_{k,l}},\frac{1}{\mu}\right)$
      for all $k$, $l$\;
    }
    $u^{(i+1)} = u^{(i,M)}$ and $v^{(i+1)} = v^{(i,M)}$\;
    $d^{u,(i+1)}_{k,l} = d^{u,(i,M)}_{k,l}$ and $d^{v,(i+1)}_{k,l} = d^{v,(i,M)}_{k,l}$
    for all $k, l$\;
    $b^{u,(i+1)}_{k,l} = b^{u,(i)}_{k,l} - d^{u,(i+1)}_{k,l} + \nabla u^{(i+1)}_{k,l}$
    for all $k, l$\;
    $b^{v,(i+1)}_{k,l} = b^{v,(i)}_{k,l} - d^{v,(i+1)}_{k,l} + \nabla v^{(i+1)}_{k,l}$
    for all $k, l$\;
  }
  \caption{The split Bregman algorithm for the OSB model.}
  \label{bh-L2L1Isoalgor}
\end{algorithm}
%
\subsection{The Model of Brox \emph{et al.}}
\label{bh-sec-broxmodel}
%
The model that we discuss in this section differs only very little from the
previous one in terms Bregman iterations. Although we have robustified the data
term and rendered the smoothness term rotationally invariant, the differences
to the previous Bregman iterative scheme will be surprisingly small. After the
discretisation of the variational formulation we obtain
\begin{equation}
  \label{eq:3}
  \argmin_{u,v} \set{\lambda D_2\left(u,v\right)  +
    \sum_{i,j} \sqrt{\norm{\nabla u_{i,j}}_2^2 + \norm{\nabla v_{i,j}}_2^2}}.
\end{equation}
Here, we observe that none of the terms of the energy functional is
differentiable. In the same way as for OSB, the smoothness term can be rewritten
in the following way
\begin{equation}
  \label{eq:15}
  \argmin_{u,v}\set{\lambda D_2\left(u,v\right) +
    \sum_{i,j} \norm{ \binom{\nabla u_{i,j}}{\nabla v_{i,j}} }_2}\enspace{}.
\end{equation}
To handle the non-differentiability, we move everything into the constraining
conditions. As before, we obtain slack variables $d^{u}$ and $d^{v}$ for the
smoothness term, and three additional slack variables $d$, $d_{x}$, and $d_{y}$
for the linearised grey value and gradient constancy assumptions in the data
term. The complete constrained optimisation task reads
\begin{align*}
  &\argmin_{\substack{u,v,\\d,d_x,d_y,\\d^u_{i,j},d^v_{i,j}}}\set{\lambda
    \left( \norm{d}_1 + \gamma \norm{\binom{d_x}{d_{y}}}_1 \right)
    + \sum_{i,j}\ \norm{\binom{d^u_{i,j}}{d^v_{i,j}}}_2} \\
  &\text{such that:}\\
  &\begin{multlined}[0.85\linewidth]
    \frac{1}{2}\left( \norm{ d - F\tbinom{u}{v} - f_t}_2^2 +
      \norm{ d_x - F_x\tbinom{u}{v} - f_{xt}}_2^2 +\right.\\
      \left.\norm{ d_y - F_y\tbinom{u}{v} - f_{yt}}_2^2 \right) +
      \frac{1}{2}\sum_{i,j}
      \norm{ \binom{d^u_{i,j} - \nabla u_{i,j}}{d^v_{i,j} - \nabla v_{i,j}} }_2^2 = 0
    \end{multlined}
\end{align*}
where the matrices $F$, $F_{x}$, and $F_{y}$ contain the spatial derivatives
from the linearised constancy assumptions. The vectors $f_{t}$, $f_{xt}$ and
$f_{yt}$ represent the corresponding spatio-temporal derivatives.
This formulation is almost identical to the one from the previous section. In
fact, the very same arguments tell us that the split Bregman algorithm is
applicable and that the minimisation with respect to $u$ and $v$ will lead to an
almost identical linear system. The biggest difference between the two
approaches lies in the minimisation with respect to $d^u$ and $d^v$. We have to
solve
\begin{equation}
  \label{eq:14}
  \argmin_{d^u,d^v}\set{\norm{\binom{d^u}{d^v}}_2 +\frac{\mu}{2}
    \norm{ \binom{d^u-\nabla u_{i,j}-b^u_{i,j}}{d^v-\nabla v_{i,j}-b^v_{i,j}} }_2^2}
\end{equation}
for each $i,j$. Fortunately, this can again be done with the help of the
generalised shrinkage operator. A detailed listing of the complete minimisation
strategy is given in Algorithm~\ref{bh-L1L1Algo}.
\begin{remark}
  Although it is possible to formulate a minimisation strategy with the Bregman
  iteration for
  $L_1$--$L_1$ models, the formulation appears a bit ``unnatural''. A few
  potential problems become immediately apparent. The first one being the fact
  that we have to eliminate the variables with respect to which we initially
  wanted to minimise completely from the cost function. Secondly, none of the
  model parameters has a direct influence on $u$ and
  $v$. They can only interact by means of the auxiliary slack variables. Chances
  are, that this will reduce the responsiveness of the algorithm to parameter
  changes. Although it is generally desirable to have algorithms that do not
  react too sensitive with respect to varying parameters, the other extreme of
  having an algorithm that reacts hardly at all, is not desirable as well.
\end{remark}
\begin{algorithm}
  \dontprintsemicolon
  \KwData{$F$, $F_x$, $F_y$, $f_t$, $f_{xt}$, $f_{yt}$, $\lambda$, $\gamma$,
    $\mu$, $N$, $M$}
  \KwResult{$u^{(N)}$ and $v^{(N)}$ minimising \eqref{eq:3}}
  \KwInit{$u^{(0)}=v^{(0)}=0$, $d^{(0)}=d_x^{(0)}=d_y^{(0)}=0$,
    ${d^{u}}^{(0)}={d^{v}}^{(0)}=0$, $f_t^{(0)}=f_t$, $=f_{xt}^{(0)}=f_{xt}$,
    $f_{yt}^{(0)}=f_{yt}$, ${b^{u}}^{,(0)}={b^{v}}^{,(0)}=0$}
  \For{$i=0$ \KwTo $N-1$}{
    $u^{(i,0)} = u^{(i)}$ and $v^{(i,0)} = v^{(i)}$ and
    $d^{(i,0)} = d^{(i)}$ and $d_x^{(i,0)} = d_x^{(i)}$ and $d_y^{(i,0)} = d_y^{(i)}$\;
    ${d^{u}}^{,(i,0)} = {d^{u}}^{,(i)}$, ${d^{v}}^{,(i,0)} = {d^{v}}^{,(i)}$\;
    \For{$j=0$ \KwTo $M-1$}{
      $\displaystyle \left(u^{(i,j+1)},v^{(i,j+1)}\right) =
      \argmin_{u,v}\left\{\norm{ d^{(i,j)} - F\tbinom{u}{v} - b^{(i)}}_2^2 + 
            \norm{ d_x^{(i,j)} - F_x\tbinom{u}{v} - b_x^{(i)}}_2^2 + \right.$\;
      $\hfill\displaystyle\frac{\mu}{2}\left( \norm{ {d^{u}}^{,(i,j)} -
              \nabla u - {b^{u}}^{,(i)} }_2^2 + \norm{{d^{v}}^{,(i,j)} -
              \nabla v - {b^{v}}^{,(i)}}_2^2\right) + $\;
      $\hfill\displaystyle\left.\norm{ d_y^{(i,j)} - F_y\tbinom{u}{v} -
              b_y^{(i)}}_2^2\right\}$\;
      $\displaystyle d^{(i,j+1)} = \shrink\left( D_x u^{(i,j+1)} + D_y v^{(i,j+1)} + 
            {f_t}^{(i)},\frac{\lambda}{\mu}\right)$\;
      $\displaystyle d_x^{(i,j+1)} = \shrink\left( D_{xx} u^{(i,j+1)} + 
            D_{xy} v^{(i,j+1)} + {f_{xt}}^{(i)},\frac{\lambda\gamma}{\mu}\right)$\;
      $\displaystyle d_y^{(i,j+1)} = \shrink\left( D_{yx} u^{(i,j+1)} +
            D_{yy} v^{(i,j+1)} + {f_{yt}}^{(i)},\frac{\lambda\gamma}{\mu}\right)$\;
      $\displaystyle \binom{d^{u,(i,j+1)}_{k,l}}{d^{v,(i,j+1)}_{k,l}} = 
      \gshrink\left( \binom{\nabla u^{(i,j+1)}_{k,l} + b^{u,i}_{k,l}}
            {\nabla v^{(i,j+1)}_{k,l} + b^{v,(i)}_{k,l}}, \frac{1}{\mu} \right)$
      for all $k$, $l$\;
    }
    $u^{(i+1)} = u^{(i,M)}$, $v^{(i+1)} = v^{(i,M)}$\;
    $d^{(i+1)} = d^{(i,M)}$, $d_x^{(i+1)} = d_x^{(i,M)}$, $d_y^{(i+1)} = d_y^{(i,M)}$\;
    ${d^{u}}^{,(i+1)} = {d^{u}}^{,(i,M)}$, ${d^{v}}^{,(i+1)} = {d^{v}}^{,(i,M)}$\;
    ${f_t}^{(i+1)} = {f_t}^{(i)} + {f_t} - d^{(i+1)} + D_x u^{(i+1)} + D_y v^{(i+1)}$\;
    ${f_{xt}}^{(i+1)} = {f_{xt}}^{(i)} + {f_{xt}} - d^{(i+1)}_x + 
    D_{xx} u^{(i+1)} + D_{xy} v^{(i+1)}$\;
    ${f_{yt}}^{(i+1)} = {f_{yt}}^{(i)} + {f_{yt}} - d^{(i+1)}_y + 
    D_{yx} u^{(i+1)} + D_{yy} v^{(i+1)} $\;
    ${b^{u}}^{,(i+1)} = {b^{u}}^{,(i)} - {d^{u}}^{,(i+1)} + \nabla u^{(i+1)}$\;
    ${b^{v}}^{,(i+1)} = {b^{v}}^{,(i)} - {d^{v}}^{,(i+1)} + \nabla v^{(i+1)}$\;
  }
  \caption{The split Bregman algorithm for the model of Brox \emph{et al.}}
  \label{bh-L1L1Algo}
\end{algorithm}
%
\section{Properties of the linear systems occurring in the Bregman algorithms}
\label{linsystsec}
%
All the linear systems that appeared in our algorithms so far have a system
matrix of the form
\begin{equation}
\label{linsystp}
\begin{gathered}
  \begin{multlined}[0.75\linewidth]
    \left(D_xD_x +  \gamma D_{xx}D_{xx} +  \gamma D_{yx}D_{yx}\right) u +\\[1mm]
    \left(D_xD_y + \gamma D_{xx}D_{xy} + \gamma D_{yx}D_{yy}\right) v - \theta \Delta u\enspace{},
  \end{multlined}\\[2mm]
  \begin{multlined}[0.75\linewidth]
    \left(D_xD_y +  \gamma D_{xx}D_{xy} +  \gamma D_{yx}D_{yy}\right) u +\\[1mm]
    \left(D_yD_y + \gamma D_{xy}D_{xy} + \gamma D_{yy}D_{yy}\right) v -
    \theta \Delta v
  \end{multlined}
\end{gathered}
\end{equation}
with parameters $\gamma > 0$, $\theta > 0$. It is interesting to note that the
discretisation of the Euler-Lagrange equations of the Horn and Schunck model
would lead to a linear system with almost the same structure. See for example
\cite{mm04}. In \cite{mm04} the authors analysed this linear system and
showed that the discretisation of the Euler-Lagrange equations leads to
symmetric and positive definite matrix. Because of the high similarity between
the two problems it will be relatively simple to adapt their proof such that we
can show the same results for our Bregman algorithms. We will even demonstrate
that the proof given in \cite{mm04} can be generalised. The authors of that
article required a specific indexing scheme for the pixels and assumed that
there was only one constancy assumption, namely the grey value constancy. The
proof given in this section demonstrates that these restrictions are not
necessary. We will show that the inclusion of higher order constancy assumptions
does not affect the positive definiteness.\par
The fact that the matrix is symmetric and positive definite is highly useful for
numerical purposes. It guarantees the convergence of algorithms such as
conjugate gradients and will allow us later on to present efficient
implementations with powerful solvers.\par
For the sake of simplicity, we will assume in the following that our image is
discretised on a rectangular grid with step sizes $h$ in each direction. We will
further assume that the pixels are indexed by a single number
$i\in\set{1,\ldots,n_p}$. The neighbouring pixels will be labelled $i_l,\ i_r,\
i_u$ and $i_d$, where the indices stand for \emph{left}, \emph{right}, \emph{up}
and \emph{down}. The sets $N_x\left(i\right)$ and $N_y\left(i\right)$ will
represent the neighbours of pixel $i$ in $x$ (resp. $y$) direction.\par
It is easy to see that the system matrix given in \eqref{linsystp} is symmetric
and positive semi-definite, since the linear system was obtained by computing
the gradient of a linear least squares system. The obtained system matrix is
large, structured and extremely sparse.\par
The first step, that we will perform, will be to rewrite the considered system
in a more explicit form. The matrices $D_x,
D_y$, etc.\ are all diagonal matrices, thus it follows that they can easily be
multiplied with each other, yielding again diagonal matrices. As for
$\Delta$, we will assume that the second derivatives are approximated in the
following way
\begin{equation}
    \left(\partial_{xx} u\right)_{i} \approx \frac{u_{i_l}-2u_{i}+u_{i_r}}{h_x^2}\enspace{},\quad
    \left(\partial_{yy} u\right)_{i} \approx \frac{u_{i_u}-2u_{i}+u_{i_d}}{h_y^2}
\end{equation}
and therefore,
\begin{equation}
    \left(\Delta u\right)_i
    = \left(\partial_{xx} u\right)_{i} + \left(\partial_{yy} u\right)_{i}
    = \sum_{k\in N_x\left(i\right)} \frac{u_{k}-u_{i}}{h^2} + \sum_{k\in N_y\left(i\right)} \frac{u_{k}-u_{i}}{h^2}\enspace{}.
\end{equation}
This leads us to the following explicit form of our linear system
($i=1,\ldots,n_p$) where the righthand side has been denoted by
$R_{u}$ and $R_{v}$ respectively.
\begin{align}
  &\begin{multlined}[0.85\linewidth]
    \left( {f_x}^2 + \gamma\left( {f_{xx}^2} + {f_{xy}^2} \right)\right)_{i} u_{i} +
    \left( {f_xf_y} + \gamma\left( {f_{xx}f_{xy}} + {f_{xy}f_{yy}}\right)\right)_{i} v_{i} -\\[1mm]
    \theta\sum_{k\in N_x\left(i\right)} \frac{u_{k}-u_{i}}{h^2}
    - \theta\sum_{k\in N_y\left(i\right)} \frac{u_{k}-u_{i}}{h^2} = R_{u_{i}}\enspace{},
  \end{multlined}\label{linsysteq1}\\[2mm]
  &\begin{multlined}[0.85\linewidth]
    \left(f_xf_y + \gamma \left(f_{xx}f_{xy} + f_{xy}f_{yy}\right)\right)_{i} u_{i} +
    \left(f_y^2 + \gamma \left(f_{xy}^2 + f^2_{yy}\right)\right)_{i} v_{i} -\\[1mm]
    \theta\sum_{k\in N_x\left(i\right)} \frac{v_{k}-v_{i}}{h^2}
    - \theta\sum_{k\in N_y\left(i\right)} \frac{v_{k}-v_{i}}{h^2} = R_{v_{i}}\enspace{}.
  \end{multlined}\label{linsysteq2}
\end{align}
In \eqref{linsystp}, the system is written down with matrices. Thus, the first
equation of \eqref{linsystp} corresponds to the $n_p$ equations given by
\eqref{linsysteq1}, whereas the second equation of \eqref{linsystp} corresponds
to the $n_p$ equations given by \eqref{linsysteq2}. If we had numbered the
equations consecutively, then \eqref{linsysteq1} would correspond to the
equations 1 to $n_p$ and \eqref{linsysteq2} would correspond to the equations
$n_p+1$ to $2n_p$. However, because of the special structure of these equations,
it is usually more convenient to write them down pairwise.\par
\begin{figure}[t]
 \centering
 \includegraphics{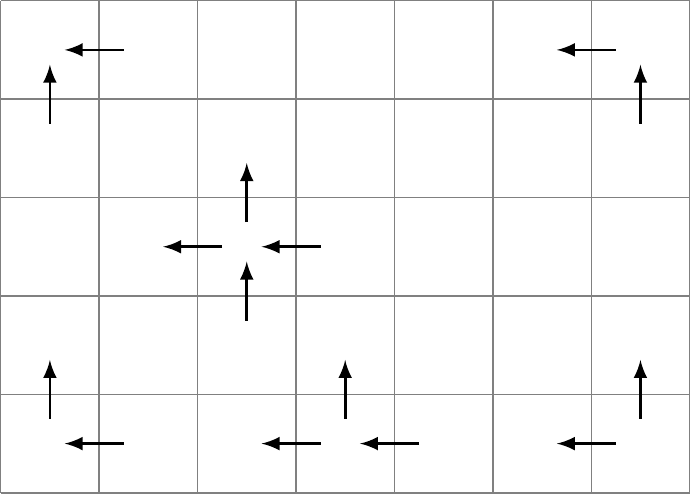}
 \caption[Reordering scheme used to show positive definiteness.]{Visualisation of the reordering scheme
 for the pixels of the image. Arrows depict to which pixels the different terms are reassigned.}
 \label{reordering}
\end{figure}

If we define the abbreviations
\begin{equation}\label{motiontensor}
  \begin{split}
    \left(J_{11}\right)_{i} &\coloneqq \left( {f_x}^2 + \gamma\left( {f_{xx}^2} + {f_{xy}^2} \right)\right)_{i}\enspace{},\\
    \left(J_{12}\right)_{i} &\coloneqq \left( {f_xf_y} + \gamma\left( {f_{xx}f_{xy}} + {f_{xy}f_{yy}}\right)\right)_{i}\enspace{},\\
    \left(J_{22}\right)_{i} &\coloneqq \left(f_y^2 + \gamma \left(f_{xy}^2 + f^2_{yy}\right)\right)_{i}
  \end{split}
\end{equation}
then we obtain the following final form
\begin{align}
  &\begin{multlined}[0.85\linewidth]
    \left( J_{11}\right)_{i} u_{i} +
    \left( J_{12} \right)_{i} v_{i} -
    \theta\sum_{k\in N_x\left(i\right)} \frac{u_{k}-u_{i}}{h^2}
    - \theta\sum_{k\in N_y\left(i\right)} \frac{u_{k}-u_{i}}{h^2} = R_{u_{i}}\enspace{},
  \end{multlined} \label{systeeq1}\\
  &\begin{multlined}[0.85\linewidth]
    \left( J_{12} \right)_{i} u_{i} +
    \left( J_{22} \right)_{i} v_{i} -
    \theta\sum_{k\in N_x\left(i\right)} \frac{v_{k}-v_{i}}{h^2}
    - \theta\sum_{k\in N_y\left(i\right)} \frac{v_{k}-v_{i}}{h^2} = R_{v_{i}}\enspace{}.
  \end{multlined}\label{systeeq2}
\end{align}
In order to show that the system matrix $M$ is positive definite, we verify that
the corresponding quadratic form $\tbinom{u}{v}^T M \tbinom{u}{v}$ is always
positive. From \eqref{systeeq1} and \eqref{systeeq2} we deduce that this
quadratic form is given by
\begin{equation}
  \begin{multlined}[0.85\linewidth]
    \sum_{i=1}^{n_p}
    \left[ \left( J_{11}\right)_{i} u_{i}^2 + 2\left( J_{12} \right)_{i}u_{i}v_{i} + 
          \left( J_{22} \right)_{i} v_{i}^2
          - \theta u_{i} \sum_{k\in N_x\left(i\right)} \frac{u_{k}-u_{i}}{h^2} \right.\\
    \left. - \theta u_{i} \sum_{k\in N_y\left(i\right)} \frac{u_{k}-u_{i}}{h^2}
          -\theta v_{i} \sum_{k\in N_x\left(i\right)} \frac{v_{k}-v_{i}}{h^2}
          -\theta v_{i} \sum_{k\in N_y\left(i\right)} \frac{v_{k}-v_{i}}{h^2}
      \right]\enspace{}.
  \end{multlined}
\end{equation}
By applying the definitions of $J_{11}$, $J_{12}$ and $J_{22}$ it is easy to see that
the first three terms in each addend can be rewritten as
\begin{equation}
  \label{linsystcond1}
  \begin{multlined}
    \left( f_{x} u + f_y v \right)^2_{i} +
    \gamma \left[ \left( f_{xx} u +f_{xy} v\right)_{i}^2 + \left( f_{xy}u+ f_{yy}v\right)^2_{i} \right]
  \end{multlined} 
\end{equation}
and thus they are always non-negative. Let us now consider the remaining terms (omitting
$\theta$ as it is strictly positive anyway):
\begin{equation}
  \label{smoothtermreorder}
  \begin{multlined}[0.85\linewidth]
    \sum_{i=1}^{n_p} \left[
          \sum_{k\in N_x\left(i\right)} \frac{u_{i}^2-u_{i}u_{k}}{h^2} +
          \sum_{k\in N_y\left(i\right)} \frac{u_{i}^2-u_{i}u_{k}}{h^2} + \right.\\
    \left.\sum_{k\in N_x\left(i\right)} \frac{v_{i}^2-v_{i}v_{k}}{h^2} +
          \sum_{k\in N_y\left(i\right)} \frac{v_{i}^2-v_{i}v_{k}}{h^2}\right]\enspace{}.
  \end{multlined}
\end{equation}
In order to show that \eqref{smoothtermreorder} is also positive, we will have
to reorder these terms one more time. This reordering is identical to the one
from \cite{mm04}. Assume that we are in pixel $i$ and that this pixel has a
neighbour in every direction. (If not, then certain terms in the following
reflection are simply not present.) Then we perform the following exchanges
(names are always based on the point of view of $i$):
\begin{itemize}
\item Pixel $i$ receives the terms $\frac{1}{h_x^2}\left( u_{i_r}^2 -u_{i}u_{i_r}\right)$
  from pixel $i_r$ and $\frac{1}{h_y^2}\left( u_{i_d}^2 -u_{i}u_{i_d}\right)$
  from pixel $i_d$.
\item Pixel $i$ receives the terms $\frac{1}{h_x^2}\left( v_{i_r}^2 -v_{i}v_{i_r}\right)$
  from pixel $i_r$ and $\frac{1}{h_y^2}\left( v_{i_d}^2 -v_{i}v_{i_d}\right)$
  from pixel $i_d$.
\item Pixel $i$ gives the terms $\frac{1}{h_x^2}\left( u_{i}^2 -u_{i}u_{i_l}\right)$
  to pixel $i_l$ and $\frac{1}{h_x^2}\left( u_{i}^2 -u_{i}u_{i_u}\right)$ to
  pixel $i_u$.
\item Pixel $i$ gives the terms $\frac{1}{h_x^2}\left( v_{i}^2 -v_{i}v_{i_l}\right)$
  to pixel $i_l$ and $\frac{1}{h_x^2}\left( v_{i}^2 -v_{i}v_{i_u}\right)$ to
  pixel $i_u$.
\end{itemize}
Figure~\ref{reordering} visualises the idea behind this reordering. The arrows
depict the direction in which a term is moved.
It follows now that \eqref{smoothtermreorder} can be rewritten as
\begin{multline}
\label{eq:10}
  \frac{1}{h^2}\sum_{i=1}^{n_p} \left[
    \smashoperator[r]{\sum_{k\in\set{i_r}\cap N_x\left(i\right)}} \left( u_i^2 + u_k^2 - 2 u_i u_k \right) \right. +
    \sum_{\mathclap{k\in\set{i_d}\cap N_y\left(i\right)}} \left( u_i^2 + u_k^2 - 2 u_i u_k \right) +\\
    \hspace*{1cm}\sum_{\mathclap{k\in\set{i_r}\cap N_x\left(i\right)}} \left( v_i^2 + v_k^2 - 2 v_i v_k \right) +
    \left.\sum_{\mathclap{k\in\set{i_d}\cap N_y\left(i\right)}}\left( v_i^2 + v_k^2 - 2 v_i v_k \right)
    \right]
\end{multline}
which is obviously always nonnegative. By applying a similar reasoning as in
\cite{mm04}, we see that \eqref{eq:10} is $0$ if and only if $u_i =
\text{const}_u$ and $v_i = \text{const}_v$ for all $i$. But then, it follows
from \eqref{linsystcond1} that $\left( f_{x} u + f_y v \right)^2_{i} = 0$ can
only be verified for all $i$ if and only if the spatial gradient $\nabla f$ is
perpendicular to the flow field $\left(u,v\right)$. In particular, it has to be
constant as well. On the other hand, $\nabla f$ is also perpendicular to the
level curves $L_c\coloneqq\set{ x\ \middle|\ f\left(x,t\right) = c }$. This
implies, that the flow field must be tangent to $L_c$ at every point. All in
all, this would mean that in the continuous setting the graph of $f$ would have
to be a plane in $\mathbb{R}^3$ for all times $t$. Thus, if we exclude this case,
where the graph is a plane, then the system matrix that we obtain in our Bregman
iterations is always positive definite.\par
The above argumentation also holds if we only consider the grey value constancy,
i.e. $\gamma = 0$. On the other hand, the grey value constancy cannot be
removed. If we only considered the constancy of the gradient, then the matrix is
not necessarily positive definite. This concludes the proof.
%
\section{Supplemental Numerical Experiments}
\label{sec:6}
\begin{figure}[tb]
  \centering
  \includegraphics[scale=0.2]{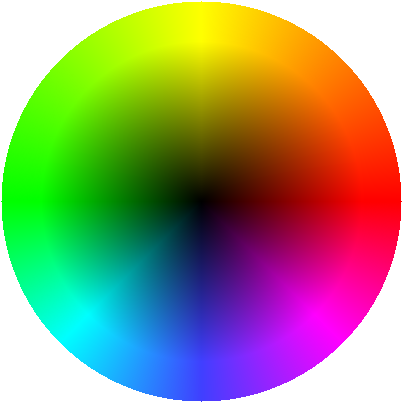}
  \caption{Color code for the displacement field.}
  \label{bh-colorcode}
\end{figure}

\begin{figure*}[t]
  \centering
  \includegraphics[width=0.213\textwidth]{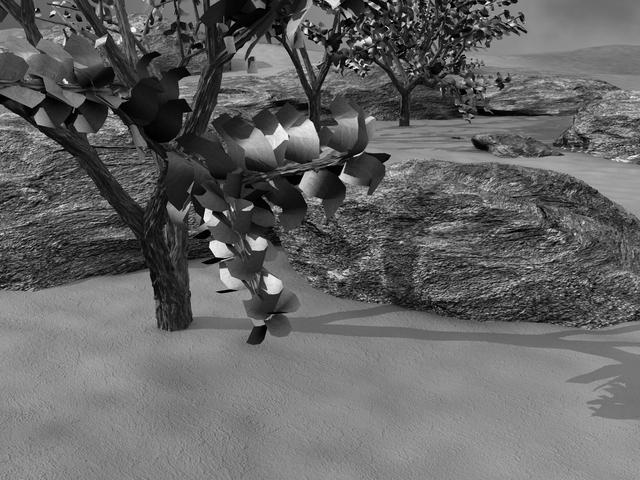}\
  \includegraphics[width=0.213\textwidth]{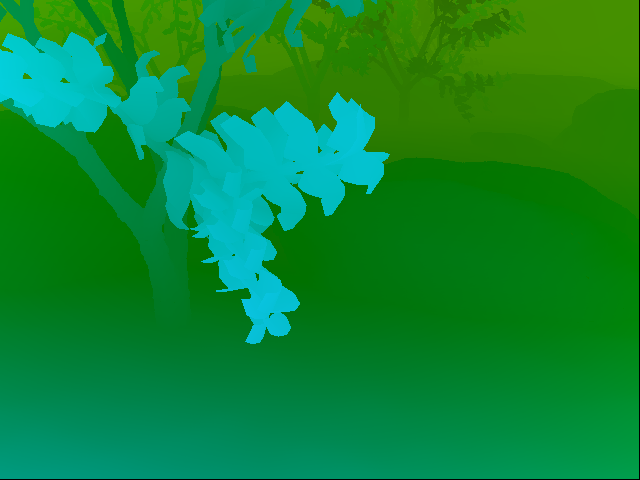}\
  \includegraphics[width=0.24\textwidth]{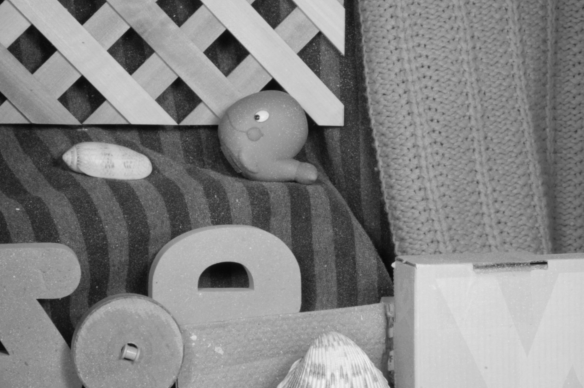}\
  \includegraphics[width=0.24\textwidth]{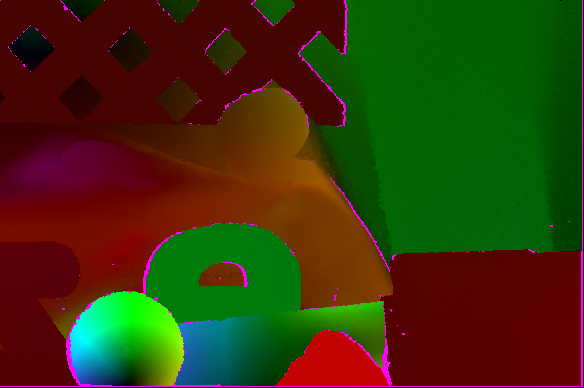}\\[1mm]
  \includegraphics[width=0.213\textwidth]{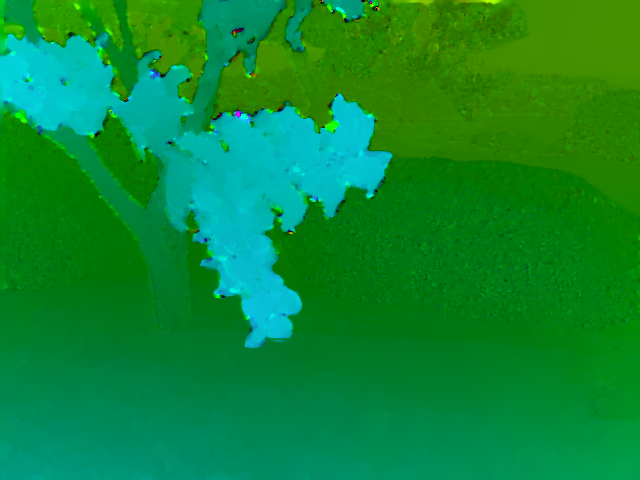}\
  \includegraphics[width=0.213\textwidth]{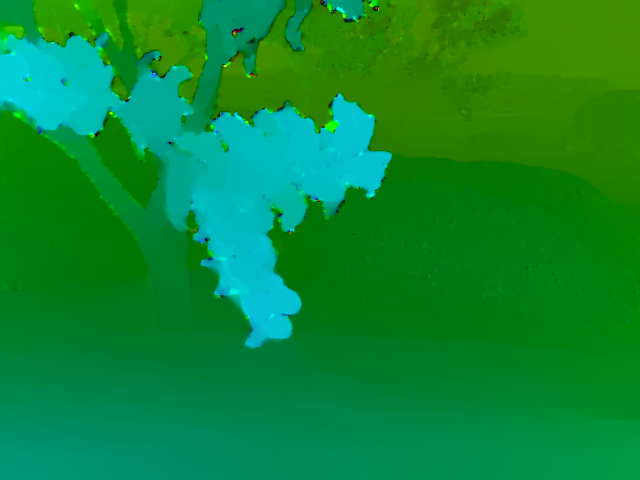}\
  \includegraphics[width=0.24\textwidth]{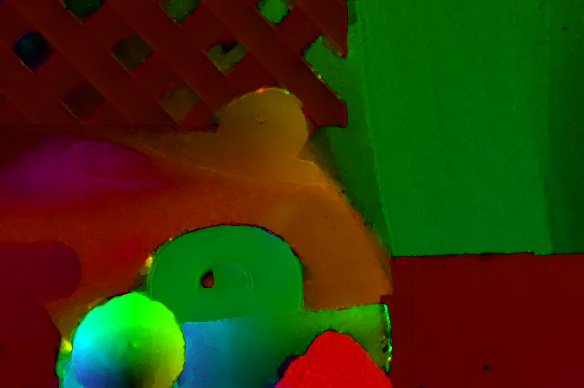}\
  \includegraphics[width=0.24\textwidth]{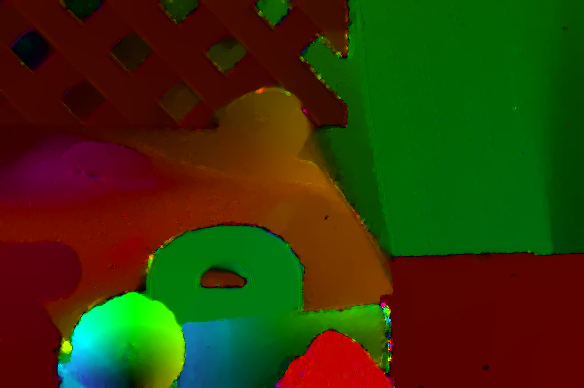}
  \caption{\emph{Top row:} The considered image sequences with their
    corresponding ground truths. \emph{From left to right:} Grove 2 frame 10,
    Grove 2 ground truth, Rubberwhale frame 10, Rubberwhale ground truth
    \emph{Bottom row:} The obtained solutions with SBM. \emph{From left to
      right:} Brox \emph{et al.}\ model, OSB model, Brox \emph{et al.}\ model,
    OSB model}
  \label{bh-image-seq}
\end{figure*}

%
The goal of this chapter will be to demonstrate the applicability of the Bregman
framework by testing some of our algorithms on a certain number of image
sequences. We will use two exemplary data sets from the Middlebury computer vision page
\cite{middlebury}. The correct ground truth of these sequences is known;
therefore, it allows us to present an accurate evaluation of our algorithms. In
order to give a quantitative representation of the accuracy of the obtained flow
fields, we will consider the so called \emph{average angular error} given by
\begin{equation}
  \text{AAE}\left(u_e,u_c\right) \coloneqq \frac{1}{n_p} \sum_{i,j}
  \arccos{\left( \frac{\langle {u_c}_{i,j} , {u_e}_{i,j} \rangle}{\norm{{u_c}_{i,j}}_2\norm{{u_e}_{i,j}}_2} \right)}
\end{equation}
as well as the \emph{average endpoint error} defined as
\begin{equation}
  \text{AEE}\left(u_e,u_c\right)\coloneqq\frac{1}{n_p}\sum_{i,j}\norm{ {u_c}_{i,j} - {u_e}_{i,j} }_2
\end{equation}
The subscripts $c$ and $e$ denote the correct respectively the estimated
spatio-temporal optic flow vectors $u_c = \left( u_{c1} , u_{c2} , 1 \right)^T$
and $u_e = \left( u_{e1} , u_{e2} , 1 \right)^T$. In this context $n_p$ denotes
the number of pixels of an image from the considered sequence. As for the
qualitative evaluation of the computed flow fields, we will use the colour
representation shown in Fig.~\ref{bh-colorcode}. Here, the hue encodes the
direction and the brightness represents the magnitude of the vector.
In the following we will use the sequences depicted in Fig.~\ref{bh-image-seq}
to test our algorithms.
The tests were done with the OSB model and the Brox \emph{et al.}\ model. The
occurring linear system was always solved with a simple Gauß-Seidel algorithm.
Usually a few dozens of iterations were more than sufficient. The number of
Bregman iterations for each model was chosen in such a way that the algorithm
reached convergence for every considered sequence.\par
Figure~\ref{bh-image-seq} depicts the obtained flow fields and the Table
\ref{bh-OSB-Results} presents the parameter choices as well as the error
measures and run times for the different algorithms. In this context RT denotes
the run time in seconds. The meanings of the parameters $\lambda$, $\mu$ and
$\gamma$ are the same as in the descriptions of the algorithms. $\sigma$ is the
standard deviation used for the preprocessing of the images with a Gaussian
convolution.\par
By looking at the results, we see that the OSB model is not only faster but also
returns the more accurate flow fields. This behaviour confirms our previous
concerns on the strong decoupling of the variables.
\begin{table}
  \centering
  \caption{Parameter choices and errors for the OSB model with 30 Bregman
    iterations, 10 Gauß-Seidel iterations and 3 alternating minimisations (first
    two rows) as well as 150 Bregman iterations, 10 Gauß-Seidel iterations and 3
    alternating minimisations (last two rows).}
  \label{bh-OSB-Results}
  \begin{tabular*}{0.95\textwidth}{@{\extracolsep{\fill}}llllllll@{}}
    \toprule
    \addlinespace
    Sequence    & $\lambda$ & $\mu$ & $\gamma$ & $\sigma$ & AAE  & AEE  & RT  \\
    \addlinespace
    \cmidrule{1-1} \cmidrule{2-5} \cmidrule{6-8} \addlinespace
    Rubberwhale & 0.0100    & 11.25 & 20.00    & 0.40     & 4.06 & 0.12 & 93  \\
    Grove 2     & 0.0250    & 6.30  & 1.50     & 0.75     & 2.79 & 0.18 & 125 \\
    \addlinespace
    Rubberwhale & 0.0065    & 0.23  & 1.00     & 0.38     & 4.67 & 0.14 & 530 \\ 
    Grove 2     & 0.0650    & 0.41  & 1.00     & 0.90     & 2.95 & 0.20 & 720 \\
    \addlinespace
    \bottomrule
  \end{tabular*}
\end{table}
%
\section{Conclusion}
\label{sec:7}
%
In this paper we have given an unified presentation of results on the Bregman
iteration that were derived from several authors in different contexts and
theoretical setups. Doing this, we have pointed out some relationships between
actual results in that field, and we have also added some new details to the
current state.\par
Furthermore, we have seen how the split Bregman algorithm can be applied to
optical flow problems. We have presented two models based on variational
formulations and showed how they can be discretised and solved with the split
Bregman algorithm. We also discussed possible algorithmical improvements like
occlusion handling and coarse-to-fine strategies that can easily be integrated
into the Bregman framework.\par
As we could see, the formulations of all the presented algorithms are very
similar. One Bregman iteration always consists in solving linear systems and
applying thresholding operations. Not only are these algorithms easy to
implement, they also offer themselves quite well for parallelisation. The
occurring linear system can be solved with Jacobi iterations, which are well
suited for massively parallel architectures such as GPUs. The shrinkage
operations are also carried out componentwise and are equally suitable for
parallel processing.\par
Finally, the positive definiteness of the system matrix allows us to consider a
broad range of highly efficient algorithms for solving the occurring linear
systems.\par
Our work is an example for a mathematical validation of important fundamental
problems in computer vision. In the future we strive to provide further
contributions in this field. Potential extensions of this work could include
anisotropic regularisers for the herein presented approaches as well as further
applications of the Bregman framework to computer vision tasks.\par
%
\bibliographystyle{abbrv}
\bibliography{bregman-arxiv}
\end{document}